\documentclass [11pt,twoside,a4paper]{article}
\usepackage{amsfonts}
\usepackage{amsthm}
\usepackage{amsmath}
\usepackage{amstext}
\usepackage{amssymb}
\usepackage{mathrsfs}
\usepackage{amscd}
\usepackage{xypic}
\usepackage{epsf}              
\usepackage{graphicx}          
\usepackage{fancybox}          
\usepackage{color}             
\usepackage{fancyhdr}
\usepackage[hang,footnotesize]{caption2}  

\def\mathcal{\mathscr}
\newfont{\aaa}{cmb10 at 19pt}
\newfont{\bbb}{cmb10 at 11pt}

\def\v1{\vspace{1mm}}

\def\leq{\leqslant}

\def\geq{\geqslant}

\newcommand{\beq}{\begin{equation}}
\newcommand{\eeq}{\end{equation}}
\newcommand{\bey}{\begin{eqnarray}}
\newcommand{\eey}{\end{eqnarray}}
\newcommand{\beyy}{\begin{eqnarray*}}
\newcommand{\eeyy}{\end{eqnarray*}}

\makeatletter
\def\@evenhead{
\vbox{\hbox to \textwidth {}{\hspace{0mm}{\footnotesize
\thepage}}{\hspace{8cm} {\footnotesize {Manli Song}}}
\protect\vspace{1truemm}\relax \hrule depth0pt
height0.15truemm width\textwidth}}
\def\@evenfoot{}
\def\@oddhead{\vbox{\hbox to \textwidth
{{\hspace{0cm}{\footnotesize Decay estimates for fractional wave equations on H-type groups}\hfill{\footnotesize
\thepage}}\hspace{0mm}}{} \protect\vspace{1truemm}\relax\hrule
depth0pt height0.15truemm width\textwidth}}
\def\@oddfoot{}
\makeatother

\begin{document}



\setcounter{page}{1}
\qquad\\[8mm]

\noindent{\aaa{Decay estimates for fractional wave equations on H-type groups}}\\[1mm]

\noindent{\bbb Manli Song}\\[-1mm]

\noindent\footnotesize{School of Natural and Applied
Sciences, Northwestern Polytechnical University, Xi'an, Shaanxi 710129, China}\\[6mm]

\normalsize\noindent{\bbb Abstract}\quad The aim of this paper is to establish the decay estimate for the fractional wave equation semigroup on H-type groups given by $e^{it\Delta^\alpha}$, $0<\alpha<1$. Combing the dispersive estimate and a standard duality argument, we also derive the corresponding Strichartz inequalties.
\vspace{0.3cm}

\footnotetext{
Corresponding author: Manli Song, E-mail:
mlsong@nwpu.edu.cn}

\noindent{\bbb Keywords}\quad Fractional wave equation, Decay estimate, H-type groups.\\
{\bbb MSC}\quad 22E25, 33C45, 35H20, 35B40\\[0.4cm]

\newtheorem{theorem}{Theorem}[section]
\newtheorem{preliminaries}{Preliminaries}[section]
\newtheorem{definition}{Difinition}[section]
\newtheorem{main result}{Main Result}[section]
\newtheorem{lemma}{Lemma}[section]
\newtheorem{proposition}{Proposition}[section]
\newtheorem{corollary}{Corollary}[section]
\newtheorem{remark}{Remark}[section]
\section{Introduction}
In this paper, we stduy the decay estimate for a class of dispersive equations:
\begin{equation}\label{fractional}
i\partial_tu+\Delta^\alpha u=f,\,u(0)=u_0,
\end{equation}
where $\Delta$ is the sub-Laplacian on H-type groups $G$, $\alpha>0$.\\

The partial differential equation in \eqref{fractional} is significantly interesting in mathematics. When $\alpha=\frac{1}{2}$, it is reduced to the wave equation; when $\alpha=1$, it is reduced to the Schr\"{o}dinger equation. The two equations are most important fundamental types of partial differential equations.

In 2000, Bahouri et al. \cite{BGX} derived the Strichartz inequalities for the wave equation on the Heisenberg group via a sharp dispersive estimate and a standard duality argument (see \cite{GV} and \cite{KT}). The dispersive estimate
\begin{equation}\label{dispersive}
||e^{it\Delta^\alpha}\varphi||_{L^\infty}\leq C|t|^{-\theta}
\end{equation}
plays a crucial role, where $\varphi$ is the kernel function on the Heisenberg group related to a Littlewood Paley decomposition introduced in Section 2 and $\theta>0$. Such estimate does not exist for the Schr\"{o}dinger equation (see \cite{BGX}). The sharp dispersive estimate is also generalized to H-type groups for the wave equation and Schr\"{o}dinger equation (see \cite{Hierro} and \cite{SZ}). Motivated by the work by Guo et al. \cite{GPW} on the Euclidean space, we consider the fractional wave equation \eqref{fractional} on H-type groups and will prove a sharp dispersive estimate.

\begin{theorem} \label{H-dispersive} Let $N$ be the homogeneous dimension of the H-type group $G$, and $p$ the dimension of its center. For $0<\alpha<1$, we have
\begin{equation*}
||e^{it\Delta^\alpha} u_0||_\infty\leq C_\alpha|t|^{-p/2}||u_0||_{\dot{B}^{N-p/2}_{1,1}},
\end{equation*}
and the result is sharp in time. Here, the constant $C_\alpha>0$ does not depend on $u_0, t$, and $\dot{B}^\rho_{q,r}$ is the homogeneous Besov space associated to the sublaplacian $\Delta$ introduced in the next section.
\end{theorem}

Following the work by Keel and Tao \cite{KT} or by Ginibre and Velo \cite{GV}, we also get a useful estimate on the solution of the fractional wave equation.
\begin{corollary}\label{Lebesgue} If $0<\alpha<1$ and $u$ is the solution of the fractional wave equation \eqref{fractional}, then for $q\in[(2N-p)/p,+\infty)$ and $r$ such that
\begin{equation*}
1/q+N/r=N/2-1,
\end{equation*}
we have the estimate
\begin{equation*}
||u||_{L^q((0,T), L^r)}\leq C_\alpha(||u_0||_{\dot{H}^1}+||f||_{L^1((0,T),\dot{H}^1)}),
\end{equation*}
where the constant $C_\alpha>0$ does not depend on $u_0$, $f$ or $T$.
\end{corollary}
\begin{remark}
In this article, we assume $0<\alpha<1$. For $\alpha=1$, the decay estimate has been proved (see \cite{Hierro}). For other cases, we could investigate the problem in the similar way to $0<\alpha\leq1$.
\end{remark}
\noindent\\[4mm]
\section{Prilimaries}
\subsection{H-type groups}
Let $\mathfrak{g}$ be a two step nilpotent Lie algebra endowed with an inner product $\langle \cdot,\cdot \rangle$. Its center is denoted by $\mathfrak{z}$. $\mathfrak{g}$ is said to be of H-type if $[\mathfrak{z}^{\bot},\mathfrak{z}^{\bot}]=\mathfrak{z}$ and for every $s \in \mathfrak{z}$, the map $J_s: \mathfrak{z}^{\bot} \rightarrow \mathfrak{z}^{\bot}$ defined by
\begin{equation*}
\langle J_s u, w \rangle:=\langle s, [u,w] \rangle, \forall u, w \in \mathfrak{z}^{\bot}.
\end{equation*}
is an orthogonal map whenever $|s|=1$.\\
\indent An H-type group is a connected and simply connected Lie group $G$ whose Lie algebra is of H-type.\\
\indent For a given $0 \neq a \in \mathfrak{z}^*$, the dual of $\mathfrak{z}$, we can define a skew-symmetric mapping $B(a)$ on $\mathfrak{z}^{\bot}$ by
\begin{equation*}
\langle B(a)u,w \rangle =a([u,w]), \forall u,w \in \mathfrak{z}^{\bot}.
\end{equation*}
We denote by $z_a$ be the element of $\mathfrak{z}$ determined by
\begin{equation*}
\langle B(a)u,w \rangle =a([u,w])=\langle J_{z_a} u,w \rangle.
\end{equation*}
Since $B(a)$ is skew symmetric and non-degenerate, the dimension of $\mathfrak{z}^{\bot}$ is even, i.e. $dim\mathfrak{z}^{\bot}=2d$.\\
For a given $0 \neq a \in \mathfrak{z}^*$, we can choose an orthonormal basis
\begin{equation*}
\{E_1(a),E_2(a),\cdots,E_d(a),\overline{E}_1(a),\overline{E}_2(a),\cdots,\overline{E}_d(a)\}
\end{equation*}
of $\mathfrak{z}^{\bot}$ such that
\begin{equation*}
B(a)E_i(a)=|z_a|J_{\frac{z_a}{|z_a|}}E_i(a)=|a|\overline{E}_i(a)
\end{equation*}
and
\begin{equation*}
B(a)\overline{E}_i(a)=-|a|E_i(a).
\end{equation*}
We set $p=dim \mathfrak{z}$. We can choose an orthonormal basis $\{\epsilon_1,\epsilon_2,\cdots,\epsilon_p \}$ of $\mathfrak{z}$ such that $a(\epsilon_1)=|a|,a(\epsilon_j)=0,j=2,3,\cdots,p$. Then we can denote the element of $\mathfrak{g}$ by
\begin{equation*}
(z,t)=(x,y,t)=\underset{i=1}{\overset{d}{\sum}}(x_i E_i+y_i \overline{E}_i )+\underset{j=1}{\overset{p}{\sum}}s_j \epsilon_j.
\end{equation*}
We identify $G$ with its Lie algebra $\mathfrak{g}$ by exponential map. The group law on H-type group $G$ has the form
\begin{equation}
(z,s)(z',s')=(z+z',s+s'+\frac{1}{2}[z,z']),  \label{equ:Law}
\end{equation}
where $[z,z']_j=\langle z,U^jz' \rangle$ for a suitable skew symmetric matrix $U^j,j=1,2,\cdots,p$.
\begin{theorem} \indent G is an H-type group with underlying manifold $\mathbb{R}^{2d+p}$, with the group law  $\eqref{equ:Law}$ and the matrix $U^j,j=1,2,\cdots$,p satisfies the following conditions:\\
$(i)$ $U^j$ is a $2d \times 2d$ skew symmetric and orthogonal matrix, $j=1,2,\cdots$,p.\\
$(ii)$ $U^i U^j+U^j U^i=0,i,j=1,2,\cdots,p$ with $i \neq j$.
\end{theorem}
{\bf Proof.} See \cite{BU}.
\begin{remark} It is well know that H-type algebras are closely related to Clifford modules (see \cite{R}). H-type algebras can be classified by the standard theory of Clifford algebras. Specially, on H-type group $G$, there is a relation between the dimension of the center and its orthogonal complement space. That is $p+1\leq 2d$ (see \cite{KR}).
\end{remark}
\begin{remark}
We identify $G$ with $\mathbb{R}^{2d}\times\mathbb{R}^p$. We shall denote the topological dimension of $G$ by $n=2d+p$. Following Folland and Stein (see \cite{FS}), we will exploit the canonical homogeneous structure, given by the family of dilations$\{\delta_r\}_{r>0}$,
\begin{equation*}
\delta_r(z,s)=(rz,r^2s).
\end{equation*}
We then define the homogeneous dimension of $G$ by $N=2d+2p$.
\end{remark}
The left invariant vector fields which agree respectively with $\frac{\partial}{\partial x_j},\frac{\partial}{\partial y_j}$ at the origin are given by
\begin{equation*}
\begin{aligned}
X_j&=\frac{\partial}{\partial x_j}+\frac{1}{2}\underset{k=1}{\overset{p}{\sum}} \left( \underset{l=1}{\overset{2d}{\sum}}z_l U_{l,j}^k \right) \frac{\partial}{\partial s_k},\\
Y_j&=\frac{\partial}{\partial y_j}+\frac{1}{2}\underset{k=1}{\overset{p}{\sum}} \left( \underset{l=1}{\overset{2d}{\sum}}z_l U_{l,j+d}^k \right) \frac{\partial}{\partial s_k},\\
\end{aligned}
\end{equation*}
where $z_l=x_l,z_{l+d}=y_l,l=1,2,\cdots,d$. In terms of these vector fields we introduce the sublaplacian $\Delta$ by
\begin{equation*}
\Delta=-\underset{j=1}{\overset{d}{\sum}}(X_j^2 +Y_j^2).
\end{equation*}
\noindent\\[2mm]
\subsection{Spherical Fourier transform}
Kor\'{a}nyi \cite{Kor}, Damek and Ricci \cite{DR} have computed the spherical functions associated to the Gelfand pair $(G, O(d))$ (we identify $O(d)$ with $O(d)\otimes Id_p$). They involve, as on the Heisenberg group, the Laguerre functions
\begin{equation*}
\mathfrak{L}_m^{(\gamma)}(\tau)=L_m^{(\gamma)}(\tau)e^{-\tau/2}, \tau \in \mathbb{R}, m,\gamma \in \mathbb{N},
\end{equation*}
where $L_m^{(\gamma)}$ is the Laguerre polynomial of type $\gamma$ and degree $m$.\\
\indent We say a function $f$ on $G$ is radial if the value of $f(z,s)$ depends only on $|z|$ and $s$. We denote respectively by $\mathcal{S}_{rad}(G)$ and $L^q_{rad}(G)$,$1\leq q\leq \infty$, the spaces of radial functions in $\mathcal{S}(G)$ and $L^p(G)$, respectively. In particular, the set of $L^1_{rad}(G)$ endowed with the convolution product
\begin{equation*}
f_1*f_2(g)=\int_Gf_1(gg'^{-1})f_2(g')\,dg', g\in G
\end{equation*}
is a commutative algebra.\\
\indent Let $f\in L^1_{rad}(G)$. We define the spherical Fourier transform, $ m\in\mathbb{N}, \lambda\in\mathbb{R}^p$,
\begin{equation*}
\hat{f}(\lambda,m)=\left( \begin{array}{c} m+d-1\\m \end{array} \right)^{-1}
\int_{\mathbb{R}^{2d+p}}e^{i\lambda s} f(z,s)\mathfrak{L}_m^{(d-1)}(\frac{|\lambda|}{2}|z|^2)\,dzds.
\end{equation*}
By a direct computation, we have $\widehat{f_1*f_2}=\hat{f_1}\cdot\hat{f_2}$. Thanks to a partial integration on the sphere $S^{p-1}$, we deduce from the Plancherel theorem on the Heisenberg group its analogue for the H-type groups.
\begin{proposition} \indent For all $f \in \mathcal{S}_{rad}(G)$ such that
\begin{equation*}
\underset{m\in\mathbb{N}}{\sum}\left( \begin{array}{c} m+d-1\\m \end{array} \right) \int_{\mathbb{R}^p} |\hat{f}(\lambda,m)||\lambda|^d d\lambda <\infty,
\end{equation*}
we have
\begin{equation}\label{plancherel}
f(z,s)=(\frac{1}{2\pi})^{d+p}\underset{m\in\mathbb{N}}{\sum}\int_{\mathbb{R}^p} e^{-i\lambda\cdot s} \hat{f}(\lambda,m) \mathfrak{L}_m^{(d-1)}(\frac{|\lambda|}{2}|z|^2)|\lambda|^d \,d\lambda,
\end{equation}
the sum being convergent in $L^{\infty}$ norm.
\end{proposition}
Moreover, if $f \in \mathcal{S}_{rad}(G)$, the functions $\Delta f$ is also in $\mathcal{S}_{rad}(G)$ and its spherical Fourier transform is given by
\begin{equation*}
\widehat{\Delta f}(\lambda,m)=(2m+d)|\lambda|\hat{f}(\lambda,m).
\end{equation*}
The sublaplacian $\Delta$ is a positive self-adjoint operator densely defined on $L^2(G)$. So by the spectral theorem, for any bounded Borel function $h$ on $\mathbb{R}$, we have
\begin{equation*}
\widehat{h(\Delta)f}(\lambda,m)=h((2m+d)|\lambda|)\hat{f}(\lambda,m).
\end{equation*}
\noindent\\[2mm]
\subsection{Homogeneous Besov spaces}
We shall recall the homogeneous Besov spaces given in \cite{Hierro}. Let $R$ be a non-negative, even function in $C_c^{\infty}(\mathbb{R})$ such that supp$R \subseteq \{\tau \in \mathbb{R}:\frac{1}{2}\leq |\tau|\leq4 \}$ and
\begin{equation*}
\underset{j\in \mathbb{Z}}{\sum}R(2^{-2j}\tau)=1, \forall \tau \neq0.
\end{equation*}

For $j\in \mathbb{Z}$, we denote by $\varphi$ and $\varphi_j$ respectively the kernel of the operator $R(\Delta)$ and $R(2^{-2j}\Delta)$. As $R \in C_0^{\infty}(\mathbb{R})$, Hulanicki \cite{Hulanicki} proved that $\varphi\in \mathcal{S}_{rad}(G)$ and
obviously $\varphi_j(z,s)=2^{Nj}\varphi(\delta_{2^j}(z,s))$.  For any $f\in\mathcal{S}^{'}(G)$, we set $\Delta_j f=f*\varphi_j$.

\indent By the spectral theorem, for any $f\in L^2(G)$, the following homogeneous Littlewood-Paley decomposition holds:
\begin{equation*}
f=\sum_{j\in\mathbb{Z}}\Delta_j f \quad \text{in $L^2(G)$}.
\end{equation*}
So
\begin{equation}\label{infty}
||f||_{L^\infty(G)}\leq\sum_{j\in\mathbb{Z}}||\Delta_j f||_{L^\infty(G)}, f\in L^2(G),
\end{equation}
where both sides of \eqref{infty} are allowed to be infinite.

Let $1\leq q,r \leq \infty, \rho <N/q$, we define the homogeneous Besov space $\dot{B}^\rho_{q,r}$ as the set of distributions $f \in \mathcal{S}^{'}(G)$ such that \begin{equation*}
||f||_{\dot{B}^\rho_{q,r}}=\left(\underset{j\in \mathbb{Z}}{\sum}2^{j\rho r}||\Delta_jf||_q^r\right)^{\frac{1}{r}}<\infty,
\end{equation*}
and $f=\underset{j\in \mathbb{Z}}{\sum}\Delta_jf$ in $\mathcal{S}^{'}(G)$.

Let $\rho <N/q$. The homogeneous Sobolev space $\dot{H}^\rho$ is
\begin{equation*}
\dot{H}^\rho=\dot{B}^0_{2,2},
\end{equation*}
which is equivalent to
\begin{equation*}
u\in\dot{H}^\rho\Leftrightarrow\Delta^{\rho/2}u\in L^2.
\end{equation*}

Analogous to Proposition 6 of \cite{FMV1} on the Heisenberg group, we list some properties of the spaces $\dot{B}^\rho_{q,r}$ in the following proposition.
\begin{proposition}\label{properties} Let $q,r\in [1,\infty]$ and $\rho<N/q$.\\
(i) The space $\dot{B}^\rho_{q,r}$ is a Banach space with the norm $||\cdot||_{\dot{B}^\rho_{q,r}}$;\\
(ii) the definition of $\dot{B}^\rho_{q,r}$ does not depend on the choice of the function $R$ in the Littlewood-Paley decomposition;\\
(iii) for $-\frac{N}{q'}<\rho<\frac{N}{q}$ the dual space of $\dot{B}^\rho_{q,r}$ is $\dot{B}^{-\rho}_{q',r'}$;\\
(iv) for any $u\in\mathcal{S}^{'}(G)$ and $\sigma>0$, then $u\in\dot{B}^\rho_{q,r}$ if and only if $L^{\sigma/2}u\in\dot{B}^{\rho-\sigma}_{q,r}$;\\
(v) for any $q_1,q_2\in[1,\infty]$, the continuous inclusion holds
\begin{equation*}
\dot{B}^{\rho_1}_{q_1,r}\subseteq\dot{B}^{\rho_2}_{q_2,r}, \quad \frac{1}{q_1}-\frac{\rho_1}{N}=\frac{1}{q_2}-\frac{\rho_2}{N}, \rho_1\geq\rho_2;
\end{equation*}
(vi) for all $q\in[2, \infty]$ we have the continuous inclusion $\dot{B}^0_{q,2}\subseteq L^q$;\\
(vii) $\dot{B}^0_{2,2}=L^2$.
\end{proposition}
\noindent\\[4mm]
\section{Technical Lemmas}
By the inversion Fourier formula \eqref{plancherel}, we may write $e^{it\Delta^\alpha}\varphi$ explicitly into a sum of a list of oscillatory integrals. In order to estimate the oscillatory integrals, we recall the stationary phase lemma.
\begin{lemma}(see \cite{S})\label{phase} Let $g \in C^\infty([a,b])$ be real-valued such that
\begin{equation*}
|g''(x)|\geq \delta
\end{equation*}
for any $x\in[a,b]$ with $\delta >0$. Then for any function $h \in C^\infty([a,b])$, there exists a constant $C$ which does not depend on $\delta, a, b, g$ or $h$, such that
\begin{equation*}
\left|\int_a^b e^{ig(x)}h(x)\,dx\right|\leq C\delta^{-1/2}\left(||h||_\infty+||h'||_1\right).
\end{equation*}
\end{lemma}

In order to prove the sharpness of the time decay in Theorem \ref{H-dispersive}. We describe the asymptotic expansion of oscillating integrals.
\begin{lemma}(see \cite{S})\label{asymptotic} Suppose $\phi$ is a smooth function on $\mathbb{R}^p$ and has a nondegenerate critical point at $\bar{\lambda}$. If $\psi$ is supported in a sufficiently small neighborhood of $\bar{\lambda}$, then
\begin{equation*}
\left|\int_{\mathbb{R}^p}e^{it\phi(\lambda)}\psi(\lambda)\,d\lambda\right| \sim |t|^{-p/2}, \text{ as t}\rightarrow \infty.\\
\end{equation*}
\end{lemma}

Besides, it will involve the Laguerre functions when we estimate the osciallatory integrals. We need the following estimates.
\begin{lemma}(see \cite{Hierro})\label{Laguerre}
\begin{equation*}
\left|(\tau \frac{d}{d\tau})^\gamma \mathfrak{L}_m^{(d-1)}(\tau)\right|\leq C_{\gamma,d}(2m+d)^{d-1/4}
\end{equation*}
for all $0\leq \gamma\leq d$.
\end{lemma}

Finally, we introduce the following properties of the Fourier transform of surface-carried measures.
\begin{theorem}(see \cite{So})\label{measure} Let $S$ be a smooth hypersurface in $\mathbb{R}^p$ with non-vanishing Gaussian curvature and $d\mu$ a $C_0^\infty$ measure on $S$. Suppose that $\Gamma\subset\mathbb{R}^p\setminus\{0\}$ is the cone consisting of all $\xi$ which are normal of some point $x\in S$ belonging to a fixed relatively compact neighborhood $\mathcal{N}$ of \text{supp} $d\mu$. Then,
\begin{align*}
\left(\frac{\partial}{\partial\xi}\right)^\nu\widehat{d\mu}(\xi)&=O\left((1+|\xi|)^{-M}\right), \forall M\in\mathbb{N},\text{ if }\xi\not\in\Gamma,\\
\widehat{d\mu}(\xi)&=\sum e^{-i(x_j,\xi)}a_j(\xi),\text{ if }\xi\in\Gamma,
\end{align*}
where the (finite) sum is taken over all points $x\in\mathcal{N}$ having $\xi$ as a normal and
\begin{equation*}
\left|\left(\frac{\partial}{\partial\xi}\right)^\nu a_j(\xi)\right|\leq C_\nu(1+|\xi|)^{-(p-1)/2-|\nu|}.
\end{equation*}
\end{theorem}

Here, we need the following properties of the Fourier transform of the measure $d\sigma$ on the sphere $S^{p-1}$. Obviously, $\widehat{d\sigma}$ is radial. By Theorem \ref{measure}, we have the radical decay properties of the Fourier transform of the spherical measure.
\begin{lemma}\label{sphere} \indent For any $\xi\in\mathbb{R}^p$, the estimate holds
\begin{equation*}
\widehat{d\sigma}(\xi)=e^{i|\xi|}\phi_+(|\xi|)+e^{-i|\xi|}\phi_-(|\xi|),
\end{equation*}
where
\begin{equation*}
|\phi^{(k)}_\pm(r)|\leq c_k(1+r)^{-(p-1)/2-k}, \text{ for all }r>0,\,k\in\mathbb{N}.
\end{equation*}
\end{lemma}
\noindent\\[4mm]
\section{Dispersive Estimates}
\begin{lemma} \label{s-t-estimate} Let $0<\alpha<1$. The kernel of $\varphi$ of $R(\Delta)$ introduced in Section 2 satisfies the estimate
\begin{equation*}
\sup_{z} |e^{it\Delta^\alpha}\varphi(z,s)|\leq C_\alpha|t|^{-1/2}|s|^{(1-p)/2}.
\end{equation*}
\end{lemma}
{\bf Proof.} By Fourier inversion \eqref{plancherel} and polar coordinate changes, we have
\begin{align}
e^{it\Delta^\alpha}\varphi(z,s)&=(\frac{1}{2\pi})^{d+p}\underset{m\in\mathbb{N}}{\sum}\int_{\mathbb{R}^p} e^{-i\lambda\cdot s+it(2m+d)^\alpha|\lambda|^\alpha}\nonumber\\
&\qquad \qquad\times R((2m+d)|\lambda|)\mathfrak{L}_m^{(d-1)}(\frac{|\lambda|}{2}|z|^2)|\lambda|^d \,d\lambda \nonumber\\
&=(\frac{1}{2\pi})^{d+p}\underset{m\in\mathbb{N}}{\sum}\int_{\mathbb{S}^{p-1}}\int_0^{+\infty} e^{-i\lambda\varepsilon\cdot s+it(2m+d)^\alpha\lambda^\alpha}\nonumber\\
&\qquad \qquad\times R((2m+d)\lambda)\mathfrak{L}_m^{(d-1)}(\frac{\lambda}{2}|z|^2)\lambda^{d+p-1} \,d\lambda d\sigma(\varepsilon).\label{after-sphere}
\end{align}
The expression after the $\mathbb{S}^{p-1}$ integral sign in \eqref{after-sphere} is very similar to an integral computed in \cite{BGX} or \cite{Hierro}( see the proof Lemma 4.1). Integrating the result over $\mathbb{S}^{p-1}$ gives us
\begin{equation}\label{general-estimate}
\sup_{z} |e^{it\Delta^\alpha}\varphi(z,s)|\leq C_\alpha\min\{1,|t|^{-1/2}\},
\end{equation}
and Lemma \ref{s-t-estimate} will come out only if we prove the case for $p\geq2$ and $|s|>1$. By switching the order of the integration in \eqref{after-sphere}, it  follows from Lemma \ref{sphere}
\begin{align*}
e^{it\Delta^\alpha}\varphi(z,s)&=(\frac{1}{2\pi})^{d+p}\underset{m\in\mathbb{N}}{\sum}\int_0^{+\infty} \widehat{d\sigma}(\lambda s)e^{it(2m+d)^\alpha\lambda^\alpha}
R((2m+d)\lambda)\\
&\qquad\qquad\qquad\qquad\times\mathfrak{L}_m^{(d-1)}(\frac{\lambda}{2}|z|^2)\lambda^{d+p-1} \,d\lambda\\
&=(\frac{1}{2\pi})^{d+p}\underset{m\in\mathbb{N}}{\sum}\int_0^{+\infty} \left(e^{i\lambda|s|}\phi_+(\lambda|s|)+e^{-i\lambda|s|}\phi_-(\lambda|s|)\right)e^{it(2m+d)^\alpha\lambda^\alpha}\\
&\qquad\qquad\qquad\qquad\times R((2m+d)\lambda)\mathfrak{L}_m^{(d-1)}(\frac{\lambda}{2}|z|^2)\lambda^{d+p-1} \,d\lambda\\
&:=(\frac{1}{2\pi})^{d+p}\underset{m\in\mathbb{N}}{\sum}\left(I_m^++I_m^-\right).
\end{align*}
Then it suffices to study
\begin{equation*}
I_m^\pm=\int_0^{+\infty} e^{i\left(\pm\lambda|s|+t(2m+d)^\alpha\lambda^\alpha\right)}\phi_\pm(\lambda|s|)R((2m+d)\lambda)\mathfrak{L}_m^{(d-1)}(\frac{\lambda}{2}|z|^2)\lambda^{d+p-1} \,d\lambda.
\end{equation*}
Performing the change of variables, $\mu=(2m+d)\lambda$, recalling that $R\in C_c^\infty(\mathbb{R})$,
\begin{equation*}
I_m^\pm=\int_{1/2}^{4} e^{itg_{m,s,t}^\pm(\mu)}h_{m,s,z}(\mu)\,d\lambda,
\end{equation*}
where
\begin{align*}
g_{m,s,t}^\pm(\mu)&=\pm\frac{\mu|s|}{(2m+d)t}+\mu^\alpha,\\
h_{m,s,z}(\mu)&=\phi_\pm\left(\frac{\mu|s|}{2m+d}\right)R(\mu)\mathfrak{L}_m^{(d-1)}\left(\frac{\mu |z|^2}{2(2m+d)}\right)\frac{\mu^{d+p-1}}{(2m+d)^{d+p}}.
\end{align*}
By Lemma \ref{Laguerre} and Lemma \ref{sphere}, we get
\begin{equation*}
||h_{m,s,z}||_\infty+||h'_{m,s,z}||_1\leq C(2m+d)^{-(2p+3)/4}|s|^{-(p-1)/2}.
\end{equation*}
Since $|\left(g_{m,s,t}^\pm\right)''|\geq \alpha|\alpha-1|2^{-\alpha-4}$, applying Lemma \ref{phase} on $I_m^\pm$ gives us
\begin{equation}\label{I-m-estimate}
\left|I_m^\pm\right|\leq C_\alpha(2m+d)^{-(2p+3)/4}|t|^{-1/2}|s|^{-(p-1)/2}.
\end{equation}
To conclude it suffices to sum these estimates since
\begin{equation*}
\underset{m\in\mathbb{N}}{\sum}(2m+d)^{-(2p+3)/4}<+\infty.
\end{equation*}

The decay estimate of time is sharp in the joint space-time cone
\begin{equation*}
\left\{(s,t)\in\mathbb{R}^p\times\mathbb{R}: s=Ct\right\}.
\end{equation*}
We will prove the sharp dispersive estimate.
\begin{lemma} \label{sharp-t-estimate} Let $0<\alpha<1$. The kernel of $\varphi$ of $R(\Delta)$ introduced in Section 2 satisfies the estimate
\begin{equation*}
\sup_{z,s} |e^{it\Delta^\alpha}\varphi(z,s)|\leq C_\alpha|t|^{-p/2}.
\end{equation*}
\end{lemma}
{\bf Proof.} From \eqref{general-estimate}, it suffices to show the inequality $|t|>1$. Recall from \eqref{after-sphere} that
\begin{equation*}
e^{it\Delta^\alpha}\varphi(z,s)=(\frac{1}{2\pi})^{d+p}\underset{m\in\mathbb{N}}{\sum}\int_{\mathbb{S}^{p-1}}I_{m,\varepsilon}\,d\sigma(\varepsilon),
\end{equation*}
where
\begin{align*}
I_{m,\varepsilon}
&=\int_0^{+\infty} e^{-i\lambda\varepsilon\cdot s+it(2m+d)^\alpha\lambda^\alpha}R((2m+d)\lambda)\mathfrak{L}_m^{(d-1)}(\frac{\lambda}{2}|z|^2)\lambda^{d+p-1} \,d\lambda\\
&=\int_{1/2}^4e^{itG_{m,\varepsilon,s,t}(\mu)}H_{m,z}(\mu)\,d\mu
\end{align*}
with
\begin{align*}
G_{m,\varepsilon,s,t}(\mu)&=\mu^\alpha-\frac{\mu}{(2m+d)t}\varepsilon\cdot s,\\
H_{m,z}(\mu)&=R(\mu)\mathfrak{L}_m^{(d-1)}\left(\frac{\mu |z|^2}{2(2m+d)}\right)\frac{\mu^{d+p-1}}{(2m+d)^{d+p}}.
\end{align*}
We will try to apply $Q$ times a non-critical phase estimate to the oscillatory integral $I_{m,\varepsilon}$.

{\bf Case 1: }$|s|\geq \alpha2^{-\alpha-3}(2m+d)|t|$. By \eqref{I-m-estimate},
\begin{equation*}
\left|\int_{\mathbb{S}^{p-1}}I_{m,\varepsilon}\,d\sigma(\varepsilon)\right|=\left|I_m^++I_m^-\right|\leq C_\alpha(2m+d)^{-p-1/4}|t|^{-p/2}.
\end{equation*}

{\bf Case 2: }$|s|\leq \alpha2^{-\alpha-3}(2m+d)|t|$. We get
\begin{equation*}
G'_{m,\varepsilon,s,t}(\mu)=\alpha\mu^{\alpha-1}-\frac{\varepsilon\cdot s}{(2m+d)t}\geq \alpha 2^{-\alpha-2}-\frac{|s|}{(2m+d)|t|}\geq \alpha2^{-\alpha-3}.
\end{equation*}
Here the phase function $G_{m,\varepsilon,s,t}$ has no critical point on $[1/2,4]$. By $Q$-fold ($1\leq Q\leq d$) integration by parts, we have
\begin{equation*}
I_{m,\varepsilon}=(it)^{-Q}\int_{1/2}^4e^{itG_{m,\varepsilon,s,t}(\mu)}D^Q\left(H_{m,z}(\mu)\right)\,d\mu
\end{equation*}
where the differential operator $D$ is defined by
\begin{equation*}
DH_{m,z}=\frac{d}{d\mu}\left(\frac{H_{m,z}(\mu)}{G'_{m,\varepsilon,s,t}(\mu)}\right).
\end{equation*}
By a direct induction,
\begin{equation*}
D^QH_{m,z}=\sum_{k=Q}^{2Q}\sum_{\wr\beta\wr=k}C(\beta, k, Q)\frac{H^{(\beta_1)}_{m,z}\left(G''_{m,\varepsilon,s,t}\right)^{\beta_2}\cdots\left(G^{(Q+1)}_{m,\varepsilon,s,t}\right)^{\beta_{Q+1}}}{\left(G'_{m,\varepsilon,s,t}\right)^k},
\end{equation*}
where $\beta=(\beta_1,\cdots,\beta_{Q+1})\in\{0,\cdots,Q\}\times\mathbb{N}^Q$ and $\wr\beta\wr=\sum_{j=1}^{Q+1}j\beta_j$.

A direct calculation shows that
\begin{equation*}
|G^{(l)}_{m,\varepsilon,s,t}(\mu)|=\alpha \prod_{j=1}^{l-1}(j-\alpha)\mu^{-l+\alpha}\leq 2^{l+2\alpha}\alpha \prod_{j=1}^{l-1}(j-\alpha)\leq C(\alpha, Q),\,l\geq2.
\end{equation*}
Using Lemma \ref{Laguerre}, by induction
\begin{equation*}
|H^{(\beta_1)}_{m,z}(\mu)|\leq C(\beta_1) (2m+d)^{-p-1/4}.
\end{equation*}
Hence, we have
\begin{equation*}
|I_{m,\varepsilon}|\leq C(\alpha, Q)|t|^{-Q}\sup_{1\leq\beta_1\leq Q}||H^{(\beta_1)}_{m,z}||_\infty\leq C(\alpha, Q)|t|^{-Q}(2m+d)^{-p-1/4}.
\end{equation*}
Taking $Q=d$, since $|t|>1$ and $p\leq 2d-1$ which implies $p/2<d$, it follows
\begin{equation*}
|I_{m,\varepsilon}|\leq C_\alpha|t|^{-p/2}(2m+d)^{-p-1/4}.
\end{equation*}
It immediately leads to
\begin{equation*}
\left|\int_{\mathbb{S}^{p-1}}I_{m,\varepsilon}\,d\sigma(\varepsilon)\right|\leq C_\alpha|t|^{-p/2}(2m+d)^{-p-1/4}.
\end{equation*}

Combining the two cases, by a straightforward summation
\begin{equation*}
|e^{it\Delta^\alpha}\varphi(z,s)|\leq C_\alpha|t|^{-p/2}\underset{m\in\mathbb{N}}{\sum}(2m+d)^{-p-1/4}\leq C_\alpha|t|^{-p/2}.
\end{equation*}
The lemma is proved.
\noindent\\[4mm]
{\bf Proof of Theorem \ref{H-dispersive}: } The dispersive inequality in Theorem \ref{H-dispersive} is a direct consequence of Lemma \ref{sharp-t-estimate} (see \cite{BGX}). It suffices to show the sharpness of the estimate.  Let $Q\in C_c^\infty\left([1/2,2]\right)$ with $Q(1)=1$. Choose $u_0$ such that
\begin{equation*}
\hat{u}_0(\lambda,m)=\begin{cases}
                      Q(|\lambda|),&m=0\\
                      0,&m\geq1.
                      \end{cases}
\end{equation*}
By the inversion Fourier formula \eqref{plancherel}, then we have
\begin{equation*}
e^{it\Delta^\alpha}u_0(z,s)=\left(\frac{1}{2\pi}\right)^{d+p}\int_{\mathbb{R}^p}e^{-i\lambda\cdot s+itd^\alpha|\lambda|^\alpha}Q(|\lambda|)e^{-|\lambda||z|^2/4}|\lambda|^d \,d\lambda.
\end{equation*}
Consider $e^{it\Delta^\alpha}u_0(0,t\bar{s})$ for a fixed $\bar{s}=\alpha d^\alpha(0,\cdots,0,1)$. The above oscillatory integral has a phase
\begin{equation*}
\Phi(\lambda)=-\lambda\cdot \bar{s}+d^\alpha|\lambda|^\alpha
\end{equation*}
with a unique nondegenerate critical point $\bar{\lambda}=\alpha^{-1} d^{-\alpha} \bar{s}=(0,\cdots,0,1)$. Indeed, the Hessian is equal to
\begin{equation*}
H(\bar{\lambda})=\alpha d^\alpha|\bar{\lambda}|^{\alpha-4}\left((\alpha-2)\bar{\lambda}_k\bar{\lambda}_l+|\bar{\lambda}|^2\delta_{k,l}\right)_{1\leq k,l\leq p}=\alpha d^\alpha\left(\begin{array}{cccc}
                           1     &    & &\\
                           &  \ddots  & &\\
                           &     &    1 &\\
                           &     &    & \alpha-1
                  \end{array}
                  \right).
\end{equation*}
So by Lemma \ref{asymptotic}, it yields
\begin{equation*}
e^{it\Delta^\alpha}u_0(0,t\bar{s})\thicksim C|t|^{-p/2}.
\end{equation*}
\noindent\\[4mm]
\section{Strichartz Inequalities}
In this section, we shall prove the Strichartz inequalities by the decay estimate in Lemma \ref{sharp-t-estimate}. We obtain the intermediate results as follows. We omit the proof and refer to \cite{GV},\cite{KT}.
\begin{theorem}\label{Intermediate} Let $0<\alpha<1$. For $i=1,2$, let $q_i, r_i\in[2,\infty]$ and $\rho_i\in\mathbb{R}$ such that
\begin{align*}
&1)\,\,2/q_i=p(1/2-1/r_i); \\
&2)\,\,\rho_i=-(N-p/2)(1/2-1/r_i),
\end{align*}
except for $(q_i, r_i, p)=(2,\infty,2)$. Let $q'_i, r'_i$ denote the conjugate exponent of $q_i, r_i$ for $i=1,2$. Then the following estimates are satisfied:
\begin{align*}
||e^{it\Delta^\alpha}u_0||_{L^{q_1}(\mathbb{R}, \dot{B}^{\rho_1}_{r_1,2})}&\leq C||u_0||_{L^2},\\
||\int_0^te^{i(t-\tau)\Delta^\alpha}f(\tau)\,d\tau||_{L^{q_1}((0,T), \dot{B}^{\rho_1}_{r_1,2})}&\leq C||f||_{L^{q'_2}((0,T),\dot{B}^{-\rho_2}_{r'_2,2})},
\end{align*}
where the constant $C>0$ does not depend on $u_0$, $f$ or $T$.
\end{theorem}

Consider the non-homogeneous fractional wave equation \eqref{fractional}. The general solution is given by
\begin{equation*}
u(t)=e^{it\Delta^\alpha}u_0-i\int_0^te^{i(t-\tau)\Delta^\alpha}f(\tau)\,d\tau.
\end{equation*}
\begin{theorem} Under the same hypotheses as in Theorem \ref{Intermediate}, the solution of the fractional wave equation \eqref{fractional} satisfies the following estimate
\begin{equation*}
||u||_{L^{q_1}((0,T), \dot{B}^{\rho_1}_{r_1,2})}\leq C(||u_0||_{L^2}+||f||_{L^{q'_2}((0,T),\dot{B}^{-\rho_2}_{r'_2,2})}),
\end{equation*}
where the constant $C>0$ does not depend on $u_0$, $f$ or $T$.
\end{theorem}

Applying Propostion \ref{properties}, by direct Besov space injections, we immediately obtain the Strichartz inequalities on Lebesgue spaces in Corollary \ref{Lebesgue}.
\noindent\\[4mm]
\noindent\bf{\footnotesize Acknowledgements}\quad\rm
{\footnotesize The work is supported by the National Natural
Science Foundation of China (Grant No. 11371036) and the Fundamental Research Funds for the Central Universities (Grant No. 3102015ZY068).}\\[4mm]

\noindent{\bbb{References}}


\begin{enumerate}
{\footnotesize
\bibitem{BGX}\label{BGX} H. Bahouri, P. G\'{e}rard et C.-J. Xu, Espaces de Besov et estimatiions de Strichartz g\'{e}n\'{e}ralis\'{e}es sur le groupe de
Heisenberg, J. Anal. Math., 2000, 82: 93-118\\[-6.5mm]
\bibitem{BU}\label{BU} A. Bonfiglioli and F. Uguzzoni, Nonlinear Liouville theorems for some critical problems on H-type groups, J. Funct. Anal., 2004, 207:
161-215\\[-6.5mm]
\bibitem{DR}\label{DR} E. Damek and F. Ricci, Harmonic analysis on solvable extensions of H-type groups, J. Geom. Anal., 1992, 2: 213-248\\[-6.5mm]
\bibitem{FS}\label{FS} G. B. Folland and E. M. Stein, Hardy spaces on homogeneous groups, Math. Notes, Princeton Univ. Press, 1992\\[-6.5mm]
\bibitem{FMV2}\label{FMV2}G. Furioli, C. Melzi and A. Veneruso, A. Littlewood-Paley decompositions and Besov spaces on Lie groups of polynomial growth,
Math. Nachr., 2006, 279: 1028-1040\\[-6.5mm]
\bibitem{FMV1}\label{FMV1}G. Furioli, C. Melzi and A. Veneruso, Strichartz inequalities for the wave equation with the full Laplacian on the Heisenberg
group,  Canad. J. Math., 2007, 59(6): 1301-1322\\[-6.5mm]
\bibitem{FV}\label{FV} G. Furioli and A. Veneruso, Strichartz inequalities for the Schr\"{o}dinger equation with the full Laplacian on the Heisenberg
group, Studia Math., 2004, 160: 157-178\\[-6.5mm]
\bibitem{GV}\label{GV} J. Ginibre and G. Velo, Generalized Strichartz inequalities for the wave equation, J. funct. Anal., 1995, 133: 50-68\\[-6.5mm]
\bibitem{GPW} Z. Guo, L. Peng and B. Wang, Decay estimates for a class of wave equations, J. Funct. Anal., 2008, 254(6): 1642–1660\\[-6.5mm]
\bibitem{Hierro}\label{Hierro} Martin Del Hierro, Dispersive and Strichartz estimates on H-type groups, Studia Math., 2005, 169: 1-20\\[-6.5mm]
\bibitem{Hulanicki}\label{Hul} A. Hulanicki, A functional calculus for Rockland operators on nilpotent Lie groups, Studia Math., 1984, 78: 253-266\\[-6.5mm]
\bibitem{Kaplan}\label{Kaplan} A. Kaplan, Fundamental solutions for a class of hypoelliptic PDE generated by composition of quadratic forms, Trans. Amer. Math.
Soc., 1980, 258: 147-153\\[-6.5mm]
\bibitem{KR}\label{HR} A. Kaplan and F. Ricci, Harmonic analysis on groups of Heisenberg type, Harmonic analysis, Lecture Nothes in Math., 1983,
992: 416-435\\[-6.5mm]
\bibitem{KT}\label{KT} M. Keel and T. Tao, Endpoints Strichartz estimates, Amer. J. Math., 1998, 120: 955-980\\[-6.5mm]
\bibitem{Kor}A. Kor\'{a}nyi, Some applications of Gelfand pairs in classical analysis, in: Harmonic analysis and group representations, 1982：333-348\\[-6.5mm]
\bibitem{LM1}\label{LM1} H. Liu and M. Song, Strichartz inequalities for the wave equation with the full Laplacian on H-type groups, Abstr. Appl. Anal., 2014, 3: 1-10\\[-6.5mm]
\bibitem{LM2}\label{LM2} H. Liu and M. Song, Strichartz inequalities for the Schr\"{o}dinger equation with the full Laplacian on H-type groups, inprint\\[-6.5mm]
\bibitem{R}\label{R} H. M. Reimann, H-type groups and Clifford modules, Adv. Appl. Clifford Algebras, 2001, 11: 277-287\\[-6.5mm]
\bibitem{So}\label{So} C. D. Sogge, Fourier integrals in classical analysis, Cambridge Univ. Press, 1993\\[-6.5mm]
\bibitem{SZ}\label{SZ} N. Song and J. Zhao, Strichartz estimates on the quaternion Heisenberg group, Bull. Sci. Math., 2014, 138(2), 293-315\\[-6.5mm]
\bibitem{S}\label{S} E. M. Stein, Harmonic analysis: real-variable methods, orthogonality and oscillatory integrals, Princeton Univ. Press, 1993\\[-6.5mm]
}
\end{enumerate}
\end{document}